\documentstyle[12pt]{article}

\begin{document}

\author{Leonid Galtchouk and Isaac Sonin \\
%EndAName
\\
IRMA, D\'epartement de Math\'ematiques,\\
Universit\'e de Strasbourg\\
7 rue R\'ene Descartes\\
Strasbourg Cedex France\\
e-mail: galtchou@math.u-strasbg.fr\\
and\\
Department of Mathematics,\\
University of North Carolina at Charlotte,\\
Charlotte, NC 28223, USA\\
e-mail: imsonin@uncc.edu \bigskip}
\title{On some estimates for bounded submartingales and the shift inequality. \\
}
\date{}
\maketitle

\begin{abstract}
It is well known that if a submartingale $X$ is bounded then the increasing
predictable process $Y$ and the martingale $M$ from the Doob decomposition $%
X=Y+M$ can be unbounded. In this paper for some classes of increasing convex
functions $f$ we will find the upper bounds for $\lim_n\sup_XEf(Y_n)$, where
the supremum is taken over all submartingales $(X_n),0\leq X_n\leq
1,n=0,1,...$. We apply the stochastic control theory to prove these results.
\end{abstract}

{\bf 1. Introduction. } Let $(\Omega ,{\cal F},({\cal F}_{n})_{n\geq 0},P)$
be a filtered probability space, $X=(X_{n})_{n\geq 0}$ be a bounded
submartingale on it, $0\leq X_{n}\leq 1,$ and let
\begin{equation}
X_{n}=Y_{n}+M_{n},n=0,1,...  \label{1}
\end{equation}
be its Doob decomposition, where $Y=(Y_{n})_{n\geq 0}$ is a predictable
nondecreasing random sequence,
\[
Y_{n}=\sum_{k=1}^{n}E(X_{k}-X_{k-1}|{\cal F}_{k-1}),\ n\geq 1,\
\]
and $M=(M_{n})_{n\geq 0}$ is a martingale. Denote by $G(0,0)$ the class of
all such submartingales with $X_{0}=Y_{0}=0,$ defined on a finite, $%
n=0,1,...,N$ or an infinite time interval, $n=0,1,...$.Though $X$ is
bounded, $Y$ and $M$ can be unbounded, respectively from above and from
below. Probably, the simplest example of such kind is when $X_{n}$ takes
only two values, $0$ and $1,$ and
\[
P(X_{n+1}=1|X_{n}=0)=P(X_{n+1}=0|X_{n}=0)=1/2,
\]
\[
P(X_{n+1}=1|X_{n})=1)=1.
\]
Then $P(Y_{n+1}-Y_{n}=1/2|X_{n}=0)=1$ and $P(Y_{n+1}-Y_{n}=0|X_{n}=1)=1.$ If
$X_{0}=Y_{0}=M_{0}=0$, then we have the Doob decomposition (\ref{1}) with
\[
Y_{n}=\frac{1}{2}+\frac{1}{2}\sum_{k=1}^{n-1}(1-X_{k}),n=1,2,....
\]

The transition probabilities defined above imply
\[
P(X_{n}=0)=(1/2)^{n},P(X_{n}=1)=1-(1/2)^{n},
\]
\[
P(Y_{n}=k/2)=P(M_{n}=-k/2)=(1/2)^{k},\ 1\leq k<n,
\]
\[
P(Y_{n}=n/2)=P(M_{n}=-n/2)=(1/2)^{n-1}.
\]
Now, one can check that $\lim_{n}Ef(Y_{n})<\infty $ if $f(x)=x^{m}$ for any $%
m\geq 0$ or if $f(x)=\exp {(\lambda x)}$ with $\lambda <2\ln 2\approx 1.386,$
and this limit is infinite for $\lambda \geq 2\ln 2.$

This example raises a natural question about the values (estimates) of
\begin{equation}
c_{n}\equiv c_{n}(f)=\sup_{X\in G(0,0)}Ef(Y_{n}),\ \hspace{0.3in}%
c(f)=\lim_{n}c_{n}(f)  \label{2}
\end{equation}
for different functions $f$. The partial answers for these questions were
given in [1], where in particular it was proved that for function $%
f(x)=x^{m} $ for any $m\geq 1$ the estimate for the upper bound is $m^{m}.$
The case of exponential $f$ was not analyzed.

The main result of this note is Theorem 1 and the inequality (\ref{5}),
which we call the {\it shift inequality.} \

{\bf Theorem 1.}

{\it (a) For functions }$f(x)=e^{\lambda x},\lambda >0$,
\[
\begin{array}{ccll}
c(f) & = & 1/(1-\lambda ), & \hbox{\rm if}\ 0<\lambda <1, \\
& = & \infty , & \hbox{\rm if}\ 1\leq \lambda .
\end{array}
\]

{\it (b) For functions }$f(x)=x^m,m\geq 1,\hspace{.3in}c(f)\leq m^m.$%
\medskip\

{\it (c) For any increasing convex function $f(x),x\geq 0,$ with concave
derivative $f^{\prime }(x)$, $c(f)\leq B<\infty $, where $B$ is the unique
solution of the equation}
\begin{equation}
B=f(0)+f^{\prime }(f^{-1}(B)).  \label{3}
\end{equation}

{\bf Remark 1. }The results from \cite{GMS} imply that $c(f)$ is finite for all
increasing convex functions with derivative of the form
\begin{equation}
f^{\prime }(x)=f^{\prime }(0)\exp \{\int_{0}^{x}\lambda (s)ds\},  \label{3b}
\end{equation}
where $\lambda (s)>0$ and $\lim_{s\rightarrow \infty }\lambda (s)=\lambda
_{0}<1,$ but the expression for $B$ is more complicated.

We will show also that the structure of submartingales $(X_{k}),$ where the
upper bounds are achieved, have the structure similar to the example above,
i.e. $X_{k}$ takes only two values, $0$ and $1,$ and $\Delta Y_{k+1}\equiv
Y_{k+1}-Y_{k}=a_{n-k},$ where the constant $a_{n-k}$ depends only on
remaining time $(n-k)$.

Before proving Theorem 1, we present briefly the main steps of the proof
grouping them into three parts A, B and C.\medskip\

{\bf A}. To obtain the estimates for $c(f)$ we consider a problem of {\it %
stochastic control} on the time interval $[0,n],$ where the control actions
are the choices of the increments $\Delta Y_{k+1}=Y_{k+1}-Y_{k}$ and $\Delta
M_{k+1}=M_{k+1}-M_{k}$ for $k=0,1,...,n-1,$ and the goal is to maximize $%
Ef(Y_{n}),$ where $f$ is a convex increasing function. Later $n$ tends to
infinity. The estimate for this functional gives an estimate for $c_{n}(f)$
and its limit gives the value of $c(f)$ in (\ref{2}).\medskip\

{\bf B}. To obtain the latter statement we show that $c(f)$ is bounded, $%
c(f)\leq lim_{n}b_{n},$ where $b_{n}$ is the solution of the following
recursive equation
\begin{equation}
b_{n}=\sup_{0\leq a\leq 1}[af(a)+(1-a)f(a+f^{-1}(b_{n-1}))].  \label{4}
\end{equation}

The sequence $b_{n}$ is nondecreasing and has a finite limit for some
functions $f$ and infinite for others. The problem here is to describe the
class of functions $f$ with finite limit. The equation (\ref{4}) is of
interest on its own, though its interpretation is not quite clear.\medskip

{\bf C}. The reduction of the problem of stochastic control to the recursive
equation (\ref{4}) is possible through the use of the following inequality
which we call the {\it shift inequality}
\begin{equation}
Ef(a+Y)\leq f(a+f^{-1}(Ef(Y))).  \label{5}
\end{equation}

Denote $S$ the class of all increasing functions $f(x),x\geq 0$ for which (%
\ref{5}) holds for any number $a\geq 0$ and any random variable $Y\geq 0$.
It is easy to see that $S$ contains all exponential functions (with equality
in (\ref{5})). We will show (Lemma 2) that this class also contains all
power functions of the form $f(x)=x^{m},m\geq 1$, and all increasing convex
functions with concave derivatives. We provide a simple example of convex
function for which (\ref{5}) is not true. We fail to find the references in
the literature to the inequalities of such type. A substantial
generalization of inequality (\ref{5}) was obtained in [2].\medskip\

We describe the problem of stochastic control in Section 2 and prove the
shift inequality for some functions $f$ in Section 3. The reduction of the
stochastic control problem to the recursive equation (\ref{4}) is given in
Section 4. Section 5 contains the solution of the recursive equation for
different classes of functions..We describe some open problems in Section 6.

We thank C. Striker who draw our attention to paper \cite{LLP} and him, F. Delbaen,
M. Emery, and J. Franchi for useful discussion concerning the shift
inequality.\medskip\

{\bf 2. A stochastic control problem}.

Similarly to $G(0,0),$ denote $G(x,y)$ = \{all submartingales $X=(X_{n})$
such that $0\leq X_{n}\leq 1,$ $X_{0}=x$, and $Y_{0}=y$, $0\leq x\leq 1,$ $%
0\leq y<\infty $\}, and given a convex increasing function $f,$%
\begin{equation}
c_{n}(f,x,y)=\sup_{X\in G(x,y)}Ef(Y_{n})=\sup_{X\in G(x,0)}Ef(y+Y_{n}),
\label{6}
\end{equation}
where $(Y_{n})$ is a predictable sequence from decomposition (\ref{1}). The
second equality in (\ref{6}) holds because $(X_{n}^{\prime })$ $\in G(x,0)$
iff $(X_{n})$ $\in G(x,y),$ where $X_{n}=X_{n}^{\prime },Y_{n}=Y_{n}^{\prime
}+y,n=0,1,...$

As in many problems of stochastic control where $n$ tends to infinity, it is
convenient to consider the maximization problem in inverse time. Thus we
consider time intervals of the form $(n,n-1,...,1,0)$ and use corresponding
notation.

Formally, we consider a Markov Decision Process (MDP) with states $(k,x,y),$
where $k=0,1,...,$ and $0\leq x\leq 1,$ $0\leq y<\infty $ (see \cite{P}). The set
of all actions (controls) {\it admissible at the state} $(k,x,y)$ is $%
\{(a,\eta ):$ $a$ is a number, $0\leq a\leq 1-x$, $\eta $ is a random
variable such that $E\eta =0$, and $0\leq x+a+\eta \leq 1\}$. We denote $%
D(x)=\{\eta :$ $E\eta =0$, $0\leq x+\eta \leq 1\}.$

The goal of control is to maximize $E_{n,x,y}f(Y_{0})$ over all possible
strategies, where a strategy is a sequence of admissible actions, maybe
nonstationary and randomized, $(n,x,y)$ is an initial state, and $Y_{0}$ is
a (random) position of a last coordinate at the last moment $k=0.$ The value
function in this problem we denote $F_{n}(x,y).$\ In notation from (\ref{2})
and (\ref{6}) $c_{n}(f)=F_{n}(0,0),$ and $c(f)=\lim_{n}F_{n}(0,0).$

Obviously
\begin{equation}
F_0(x,y)=F_n(1,y)=f(y),\hspace{.1in}n=0,1,...  \label{7}
\end{equation}

The Bellman equation for our problem takes a form
\begin{equation}
F_{n}(x,y)=\sup_{0\leq a\leq 1-x,}\{\sup_{\eta \in D(x+a)}EF_{n-1}(x+a+\eta
,y+a)\}\equiv sup_{a}T^{a}F_{n-1}(x,y),  \label{8}
\end{equation}
\medskip\ where operators $\ T^{a}g(x,y)=sup_{\eta }T^{a,\eta }g(x,y),$ $%
T^{a,\eta }g(x,y)=Eg(x+a+\eta ,y+a).$\medskip

Before to describe the structure of optimal actions in Lemma 1,\ we need the
following simple statement. Let us denote $\xi (x)\in D(x)$ a random
variable taking value $(1-x)$\ with probability $x,$\ and value $-x$\ with
probability $(1-x).$\ \medskip\ \

{\bf Proposition 1.} {\it Let $h(s)$ be a convex function on $[0,1],x\in
[0,1].$ Then }

{\it $\sup_{\eta \in D(x)}Eh(x+\eta )=Eh(x+\xi (x))\equiv xh(1)+(1-x)h(0). $}%
\medskip

Proof. The set $D(x)$ is convex and closed in a weak topology. The random
variable $\xi (x)$ is an unique extreme point and hence $Ef(x+\eta )$
reaches its maximum at $\xi (x).$ \medskip\

{\bf Lemma 1. }{\it For the problem of stochastic control described above on
a finite or an infinite time interval}

{\it \ a) there is an optimal action $(a_n(x,y),\eta _n(x,y))$ at each state
$(n,x,y)$, where }$n=0,1,2,...${\it \ is a remaining time,}

{\it b) the second component of the optimal action has a form }$\eta ${\it $%
_{n}(x,y)=\xi (x+a_{n}(x,y))$, where the random variable $\xi (u)$ is
described above.}

{\it c) $F_{n}(x,y)$ is an increasing function in }$y${\it \ given }$n,x,$%
{\it \ and a decreasing and convex in $x$ (constant for $n=0$) given $n,y.$}%
\medskip\

Proof. Point a) follows from the general theory of MDP. Point b) says that $%
X_{n}=0$ or $1$ for all $n$, except maybe the initial moment, with
probability $1$. We will prove points b) and c) by induction on $n.$

For $n=0$ we have $F_{0}(x,y)=f(y)$, i.e. c) holds. Suppose that Lemma 1 is
proved for $n-1.$ Since $F_{n-1}(x,y)$ is convex in $x$, Proposition 1
immediately implies point b) of Lemma 1 for $n$. Then, using also the second
equality in (\ref{7}) for $(n-1),$ we obtain that in the Bellman equation (%
\ref{8})
\begin{equation}
T^{a}F_{n-1}(x,y)=[(x+a)f(y+a)+(1-(x+a))F_{n-1}(0,y+a)].  \label{9a}
\end{equation}
and hence the Bellman equation takes a form
\begin{equation}
F_{n}(x,y)=\sup_{0\leq a\leq 1-x}[(x+a)f(y+a)+(1-(x+a))F_{n-1}(0,y+a)].
\label{9}
\end{equation}

Since $F_{n-1}(0,y)$ and $f(y)$ are increasing in $y,$ formula (\ref{9})\
implies that $F_{n}(x,y)$ is also increasing in $y.$

To prove the monotonicity property of $F_{n}(x,y)$ in $x,$ let us note that
if $0\leq x_{1}<x_{2}\leq 1,$ and action $(a,\eta )$ is admissible at $x_{2}$
then action $(a+x_{2}-x_{1},\eta )$ is admissible at $x_{1}$. Then formula (%
\ref{9a}) implies that $%
T^{a+x_{2}-x_{1}}F_{n-1}(x_{1},y)=T^{a}F_{n-1}(x_{2},y)+(x_{2}+a)[f(y+a+x_{2}-x_{1})-f(y+a)]+(1-(x_{2}+a))[F_{n-1}(0,y+a+x_{2}-x_{1})-F_{n-1}(0,y+a)].
$ The monotonicity of functions $F_{n-1}(x,y)$ and $f(y)$ in $y$ implies
that the latter expression is positive. Since this is true for any $a$
admissible at $x_{2},$ we obtain that $F_{n}(x_{1},y)>F_{n}(x_{2},y).$

Now let us prove the convexity of $F_{n}(x,y)$ in $x$, i.e. the inequality
\begin{equation}
(F_{n}(x_{1},y)+F_{n}(x_{2},y))/2\geq F_{n}(x,y),\hspace{0.1in}%
x_{1}+x_{2}=2x,0\leq x_{1}<x_{2}\leq 1.  \label{10}
\end{equation}
Let an action $a_{n}(x,y)\equiv a$ is an optimal action at state $(n,x,y)$,
and thus $0\leq a\leq 1-x.$

Suppose first that $a=1-x$. Then by Lemma 1 $\eta _{n}(x,y)=\xi (1)\equiv 0$
and by (\ref{9}) and the second equality in (\ref{7})
\begin{equation}
F_{n}(x,y)=F_{n-1}(1,y+1-x)=f(y+a).  \label{11}
\end{equation}
Let us consider any points $x_{1},x_{2}$ such that $x_{1}+x_{2}=2x$. Denote $%
a_{i}=1-x_{i}$. Then $a=1-x=(a_{1}+a_{2})/2$. Using formula (\ref{9a}) for $%
x=x_{i},i=1,2$ we obtain\medskip\

$F_n(x_i,y)\geq T^{a_i}F_{n-1}(x_i,y)=F_{n-1}(1,y+a_i)=f(y+a_i),i=1,2$%
.\medskip\

By convexity of function $f$ we have $(f(y+a_{1})+f(y+a_{2}))/2\geq f(y+a)$,
and therefore, using (\ref{11}), we obtain (\ref{10}).

Suppose now that $a<1-x$. In this case we prove (\ref{10}) for $x_{1},x_{2},$
such that $x_{1}+x_{2}=2x$ and $|x_{i}-x|/2\leq 1-(a+x),i=1,2.$. This of
course implies the convexity of $F_{n}(x)$. Note that for such $x_{i}$ an
action $(a,\xi (x_{i}+a))$ is an admissible action. Let us show that for
such $x_{1},x_{2}$
\begin{equation}
\lbrack T^{a,\xi (x_{1}+a)}F_{n-1}(x_{1},y)+T^{a,\xi
(x_{2}+a)}F_{n-1}(x_{2},y)]/2=F_{n}(x,y)=T^{a,\xi (x+a)}F_{n-1}(x,y).
\label{12}
\end{equation}

By the definition of operators $T^{a,\eta }$ we have
\[
T^{a,\xi (x+a)}F_{n-1}(x,y)=(x+a)f(y+a)+(1-(x+a))F_{n-1}(0,y+a),
\]
and

$T^{a,\xi
(x_i+a)}F_{n-1}(x_i,y)=(x_i+a)f(y+a)+(1-(x_i+a))F_{n-1}(0,y+a),i=1,2.$%
\medskip\

Taking the average of the last two equalities and using the equality $%
x_{1}+x_{2}=2x$ we obtain (\ref{12}). Since $F_{n}(x_{i},y)\geq T^{a,\eta
}F_{n-1}(x_{i},y)$ for all admissible $(a,\eta )$, then (\ref{12}) implies (%
\ref{10}). Lemma 1 is proved.\medskip\

{\bf 3. The Shift inequality (\ref{5}). The description of the class }$S.$%
\medskip\

Now we turn our attention to the description of the class $S$ of all
increasing functions $f(x),x\geq 0,$ not necessarily convex, for which (\ref
{5}) holds.\medskip\ In \cite{GMS} it was proved that a necessary and sufficient
condition for twice continuously differentiable functions $f$ to be in $S$ :
$\frac{f^{\prime \prime }(x)}{f^{\prime }(x)}$ is a nonincreasing function.
\

To keep our paper selfcontained we present a brief and different proof for
the cases covered in our Theorem 1.\smallskip

{\bf Lemma 2. }{\it Class \ }$S${\it \ contains }

a) {\it all exponential functions }$f(x)=e^{\lambda x},\lambda >0,$

b) {\it all power functions}$\ f(x)=x^m,m\geq 1$, {\it \ and}

c){\it \ all increasing convex functions with concave derivatives.}\medskip\

Proof. The first statement is checked trivially.

To prove b) and c) note that we can rewrite (\ref{5}) as

\begin{equation}
f^{-1}(Ef(a+Y))\leq a+f^{-1}(Ef(Y))  \label{13}
\end{equation}
and one has the equality if $a=0$.

b) Let $f(x)=x^{m},m\geq 1,$ and $\Vert X\Vert _{m}=(E|X|^{m})^{1/m}.$ Then
by the triangle (Minkovski) inequality
\[
f^{-1}(Ef(a+Y))=\Vert a+Y\Vert _{m}\leq a+\Vert Y\Vert _{m}=a+f^{-1}(Ef(Y)).
\]

c)\ The derivative (with respect to $a$) of the left-hand side term in (\ref
{13}) is equal to
\[
\lbrack f^{-1}(Ef(a+Y)]^{\prime }=\frac{Ef^{\prime }(a+Y)}{f^{\prime
}(f^{-1}(Ef(a+Y)))},
\]
where in the numerator we use $(Ef(a+Y))^{\prime }=Ef^{\prime }(a+Y)$ since $%
f$ is continuously differentiable. Since $f$ is a convex and increasing
function, we have $Ef(a+Y)\geq f(E(a+Y)),$ and therefore $%
f^{-1}(Ef(a+Y))\geq E(a+Y),$ and $f^{\prime }(f^{-1}(Ef(a+Y)))\geq $ $%
f^{\prime }(E(a+Y))\geq Ef^{\prime \prime }(a+Y),$ where the last inequality
is true since $f^{\prime }$ is concave. Thus we obtained that $%
[f^{-1}(Ef(a+Y)]^{\prime }\leq 1$ and therefore (\ref{13}). Lemma 2 is
proved.\medskip\

{\bf Remark 2. }Two last cases may suggest a conjecture that the shift
inequality would be true for all convex slowly increasing functions. It is
easy to see that the inequality is wrong for the following slowly increasing
function : $f(x)=x$ if $0\leq x\leq 1$,\ $f(x)=(1+x^{2})/2$ if $1\leq x\leq
\infty ,$ and $a=1,\ Y$\ is a symmetrical Bernoulli random
variable.\medskip\

{\bf 4. From the Bellman equation (\ref{8}) to the recursive equation (\ref
{4}) through the shift inequality (\ref{5}). }\medskip\

Our final goal is to estimate $c(f)=\lim_nF_n(0,0).$ By Lemma 1 and (\ref{9}%
) we have
\begin{equation}
F_n(0,0)=\sup_{0\leq a\leq 1}(af(a)+(1-a)(F_{n-1}(0,a)).  \label{14}
\end{equation}

By definition $F_{n-1}(0,0)=\sup_{\pi }Ef(Y_{0}),$ where $\pi $ is an
admissible strategy and $Y_{0}$ is a corresponding position of $y$. Denote $%
H(n-1,0,0)$ the set of all possible $Y_{0}$ for an initial point $(n-1,0,0).$
The set of all admissible strategies for any initial point $(k,x,y)$ does
not depend on $y$. Hence $H(n-1,0,a)=a+H(n-1,0,0)$ and $F_{n-1}(0,a)=%
\sup_{H(n-1,0,a)}Ef(Y_{0})=\sup_{H(n-1,0,0)}Ef(a+Y_{0}).$ This statement is
just a paraphrase of a second equality in (\ref{6}).

Suppose that a convex function $f$ belongs to the class $S,$ for which the
shift inequality (\ref{5}) holds. Then $Ef(y+Y_{0})\leq
f(y+f^{-1}(Ef(Y_{0})))$ and therefore, using the fact that both functions $f$
and $f^{-1}$ are strictly increasing, we have
\begin{eqnarray*}
F_{n-1}(0,a) &=&\sup_{\pi }Ef(a+Y_{0})\leq \sup_{\pi }f(a+f^{-1}(Ef(Y_{0}))=
\\
&=&f(a+f^{-1}(\sup_{\pi }Ef(Y_{0}))=f(a+f^{-1}(F_{n-1}(0,0)).
\end{eqnarray*}

Combining this with (\ref{14}) and using notation $F_n(0,0)=c_n,$ we obtain
the inequality
\begin{equation}
c_n\leq \sup_{0\leq a\leq 1}[af(a)+(1-a)f(a+f^{-1}(c_{n-1}))].  \label{15}
\end{equation}

Comparing this sequence with a sequence $(b_{n})$ defined by the recursive
equation (\ref{4}), and assuming that $b_{0}=c_{0}=f(0),$ we obtain that $%
c_{n}\leq b_{n}$ for all $n,$ and hence the upper estimate for $\lim b_{n}$
can serve as an estimate for $\lim c_{n}$. Note that for exponential
functions inequality in (\ref{15}) become an equality and therefore $%
c_{n}=b_{n}$ for all $n$\ and $c(f)=\lim b_{n}.$\medskip\

{\bf 5. The solution of the recursive equation (\ref{4}). } \medskip\

Given a convex increasing function $f,$ let us consider the function $%
g(a,b)=af(a)+(1-a)f(a+f^{-1}(b))$ and the function $G(b)=\sup_{0\leq a\leq
1}g(a,b).$ Then $b_{n}=G(b_{n-1}),n=1,2,...,b_{0}=f(0)$. Since $g(0,b)=b$ we
have $b\leq G(b)$ and $b_{n-1}\leq b_{n}.$

Denote
\begin{equation}
B=\inf \{b\geq 0:\partial g(a,b)/\partial a\leq 0\ \hbox{\rm for all}\ 0\leq
a\leq 1\}.  \label{16}
\end{equation}
Since \ $\partial g(a,B)/\partial a$ ${\leq 0}$ for all $a,$ and $g(0,B){=B}$
we have $\sup_ag(a,B)=B.$ Function $g(a,b)$ is increasing in $b$ for any
fixed $a.$ Hence, if $b<B$ then $\sup_ag(a,b)\leq \sup_ag(a,B)=B.$ So we
immediately obtain the following \medskip\

{\bf Proposition 2.} {\it A sequence }$b_n${\it \ is increasing and }$%
\lim_nb_n\leq B.$\medskip\

Our goal is to describe functions $f$ for which $B<\infty $ and to find
conditions when $\lim b_n=B.$ We have
\begin{equation}
\partial g(a,b)/\partial a \, {\equiv g^{\prime }(a,b)=f(a)+af^{\prime
}(a)-f(a+f^{-1}(b))+(1-a)}f^{\prime }(a+f^{-1}(b)).  \label{17}
\end{equation}

In particular {$g^{\prime }(0,b)=f(0)-b$ }${+}$ $f^{\prime }(f^{-1}(b)).$%
\medskip\

{\bf Proposition 3.} {\it If }$g^{\prime }(a,b)\geq c>0,${\it \ for all }$%
0\leq a<${\it \ }$\varepsilon ${\it , }$\varepsilon >0,\ ${\it then }$%
G(b)\geq b+c\varepsilon $.\medskip\

Proof. The definition of $G(b)$ and a condition of Proposition 3 imply that $%
G(b)\geq g(\varepsilon ,b)\geq g(0,b)+c\varepsilon =b+c\varepsilon .$%
\medskip\

Proposition 3 implies immediately a simple sufficient condition for $%
\lim_nb_n=\infty $.\medskip\

{\bf Proposition 4.} {\it If for some }$\varepsilon >0,${\it \ }$g^{\prime
}(a,b)\geq c>0,${\it \ for all }$0\leq a<\varepsilon ${\it \ and all }$b>0,$%
{\it \ then }$\lim_{n}b_{n}=\infty .$\medskip\

Proof of Theorem 1. (a) For function $f(x)=e^{\lambda x},\lambda >0$, we
have
\[
g^{\prime }(a,b)=[e^{\lambda a}(a+(1-a)b)]^{\prime }=e^{\lambda
a}(1+b(\lambda (1-a)-1)+\lambda a)
\]
and
\[
g^{\prime }(0,b)=1+b(\lambda -1).
\]
Therefore:

if $\lambda >1$ then it is easy to see that the condition of Proposition 4$%
\; $holds for sufficiently small $\varepsilon ,$ and we obtain that $%
\lim_nb_n=\infty .$

if $\lambda <1$ then it is easy to see that $B=1/(1-\lambda )$ and if $%
\lim_nb_n=b_{*}<B$ we obtain a contradiction with Proposition 3.

If $\lambda =1$ then $g^{\prime }(a,b)=e^{a}(1+a(1-b))$ and starting with $%
b_{0}=f(0)$ we obtain $b_{1}=G(b_{0})=f(1)=e.$ For $b\geq 2$ we have $%
\sup_{a}g(a,b)=g(1/(b-1),b)=e^{1/(b-1)}(b-1)$ and if $\lim_{n}b_{n}=b_{\ast
}<\infty $ then $b_{\ast }$ must satisfy the equation $b=(b-1)e^{1/(b-1)}.$
It is easy to check that this equation has no solution. Thus $%
\lim_{n}b_{n}=\infty $.

(b) Consider the function $f(x)=x^m,m\geq 1.$ Formula (\ref{17}) gives
\[
g^{\prime }(a,b)=(m+1)(a^m-a(a+b^{1/m})^{m-1})+(a+b^{1/m})^{m-1}(m-b^{1/m}).
\]
Therefore $B=m^m$ and $\lim_nb_n=b_{*}\leq B$. It is easy to check that if $%
b<B$ then there are $c>0,\varepsilon >0$ such that $g^{\prime }(a,b)\geq c>0$
for all $0\leq a<\varepsilon $ which contradicts to Proposition 3. Hence $%
b_{*}=B.$

(c) Formula (\ref{17}) can be rewritten as
\begin{eqnarray}
g^{\prime }(a,b) &=&f(a)+a[f^{\prime }(a)-f^{\prime
}(a+f^{-1}(b))]+[f^{\prime }(a+f^{-1}(b))-f(a+f^{-1}(b))]  \label{18} \\
&=&f(a)+ah_1(a,b)+h_2(a,b).  \nonumber
\end{eqnarray}
Since $f^{\prime }$ is increasing we have $h_1(a,b)<0$. Since $f$ is a
convex and increasing and $f^{\prime }$ is a concave then $h_2(a,b)<0$ for
sufficiently large $b$, and tends to $-\infty $ as $b\to \infty $. It
implies that for sufficiently large $b$ the derivative $g^{\prime }(a,b)<0$
for all $0\le a\le 1$. Hence $B<\infty ,\ $where $B$ is defined in (\ref{16}%
). As in point (b) one shows that the inequality $\lim_nb_n=b_{*}<B$ implies
a contradiction with Proposition 3. Since $h_2^{\prime }(a,b)<0$ for all $a,$
we obtain also that the value of $B$ defined in (\ref{16}) can be found from
the condition $g^{\prime }(0,b)=0.$ In other words the value of $B$ defined
in (\ref{16}) is also a solution of an equation $B=f(0)+f^{\prime
}(f^{-1}(B)).$\medskip\ \

{\bf 6. Some solved and some open problems.}\medskip

The following problem was solved in \cite{GMS}.

{\bf Problem 1.} Describe the class of all increasing functions $f(x),x\geq
0,$ not necessarily convex, for which (\ref{5}) holds.

The following problems are open.

{\bf Problem 2.} Find precise estimates for $c(f)$ in cases b) and c) of
Theorem 1.

{\bf Problem 3.} Obtain results about the possible growth of $Ef(Y_{n})$
when instead of boundedness of $(X_{n})$ some assumptions on its growth are
imposed.\

{\bf Problem 4.} Find an independent and natural interpretation of the
recursive equation (\ref{4}).


\begin{thebibliography}{99}
\bibitem{LLP}   Lenglart, E., Lepingle, D., and Pratelli, M. (1980). Pr\'{e}sentation
unifi\'{e}e de certaines in\'{e}galit\'{e}s de la th\'{e}orie des
martingales. {\sl Seminaire de Probabilit\'{e}s. Lect. Notes in Math.} {\bf %
784}, 26-48.

\bibitem{GMS} Gordon, A., Molchanov, S., and Sonin, I. (2010). The shift and index inequalities
and their properties, submitted to JIPAM.

\bibitem{P} Puterman, M. (1994). {\sl Markov Decision processes.} Wiley \& Sons,
New-York.

\end{thebibliography}
\end{document}